\documentstyle[12pt]{article}

\newcommand{\ts}{\otimes} 

\newcommand{\pp}{{\bf P}}
\newcommand{\qq}{{\bf Q}}
\newcommand{\rr}{{\bf R}}
\newcommand{\zz}{{\bf Z}}
\newcommand{\cc}{{\bf C}}
\newcommand{\oo}{{\cal O}}
\newcommand{\ii}{{\cal I}}

\newcommand{\ra}{\rightarrow}

\newcommand{\ord}{{\rm ord}}

\newcommand{\base}{{\rm BS}}
\newcommand{\mult}{{\rm mult}}

\newtheorem{theorem}{Theorem}[section]
\newtheorem{lemma}[theorem]{Lemma}

\newtheorem{conjecture}[theorem]{Conjecture}
\newtheorem{corollary}[theorem]{Corollary}
\newtheorem{proposition}[theorem]{Proposition}
\newtheorem{definition}[theorem]{Definition}
\newcounter{nmb}

\newcommand{\next}{\addtocounter{nmb}{1}}
\newcommand{\rn}{\addtocounter{theorem}{1}}

\newcommand{\nonmbthm}{\renewcommand{\theequation}{\thetheorem}}

\begin{document}
\setcounter{section}{-1}
\title{Seshadri constants at very general points}
\author{Michael Nakamaye}
\maketitle

\section{Introduction}
\nonmbthm
The goal of this paper is to study the Seshadri constant of an ample line
bundle $A$ at a very general point $\eta$ of a smooth projective variety $X$.  
\begin{definition} Suppose $X$ is a smooth projective variety, $x \in X$, and
$A$ an ample line bundle on $X$. Then we define the Seshadri constant of 
$A$ at $x$ by
$$
\epsilon(x,A) = \inf_{C \ni x}\frac{c_1(A) \cap C}{\mult_x(C)};
$$
here the infimum runs over all integral
curves $C \subset X$ passing through $x$. 
\label{d1}
\end{definition}
Equivalently, if $\pi: Y \ra X$ denotes the blow--up of $X$ at $x$ with
execptional divisor $E$, then
$$
\epsilon(x,A) = \sup_r \{r \in \qq^+ \,|\, \pi^\ast(A)(-rE) \,\, \mbox{is ample}\}.
$$
The Seshadri constant $\epsilon(x,A)$ measures how many jets $nA$ 
separates at $\eta$ asymptotically as $n \ra \infty$.  
In the case when $X$ is a surface, it is known \cite{el} that
$\epsilon(\eta,A) \geq 1$.  
Meanwhile,
Ein, K\"uchle, and Lazarsfeld \cite{ekl} have established a lower bound
in arbitrary dimension
$$
\epsilon(\eta,A) \geq \frac{1}{\dim X}.
$$
The factor of $\dim X$ appearing in the general
result is a function of the ``gap argument'' used in the proof.
The same gap argument is also responsible for the presumably extra
factor of $\dim X$ in the known results for global generation of adjoint 
bundles (see \cite{S} for example).  

Our goal in this paper is to make some progress toward obtaining 
better lower bounds for $\epsilon(\eta,A)$.  Ultimately, however,
one would like a bound which does not
depend on the dimension of $X$ and so from our point of view the interest of
this work lies more in the methods used than in the explicit results.  
The basic idea employed here, 
namely that a singular Seshadri exceptional subvariety influences the
dimension count for sections with specified jets, 
was already presented in \cite{N}.  The counting is difficult and we have 
tried to find a compromise between computational complexity versus obtaining
the best possible results.  Thus we have counted very carefully in the 
three--fold case and less so in the higher dimensional case.

\medskip

In order to clarify our basic strategy we will
go through the argument completely 
in the three--fold case to obtain:  

\begin{theorem}
Suppose $X$ is a smooth three--fold and $\eta \in X$ a very general point.
Then for any ample line bundle $A$ on $X$ we have
$$
\epsilon(\eta,A) \geq \frac{1}{2}.
$$
\label{t}
\end{theorem}

The proof of Theorem \ref{t} will use the work of
Ein and Lazarsfeld \cite{el} on surfaces 
as well as the uniform bounds on symbolic
powers obtained in \cite{els}.  
We will then establish the following slight quantitative improvement of the
main result in \cite{ekl} for arbitrary dimension:

\begin{theorem}  Suppose $X$ is a projective variety of dimension $d \geq 4$ 
and
$A$ an ample line bundle on $X$.  Then for a very general point $\eta \in X$
we have the lower bound
$$
\epsilon(\eta,A) > \frac{3d+1}{3d^2}.
$$
\label{main}
\end{theorem}

\vspace{-.15cm}
The proofs of Theorem \ref{t} and Theorem \ref{main} follow \cite{ekl} 
line for line, our only innovation coming in counting jets.  The fundamental
observation is that if $\epsilon(\eta,A)$ is small this puts restrictions
on $h^0\left(X,nA \ts m_\eta^{n\alpha}/m_\eta^{n\alpha+1}\right)$ 
for various values of $\alpha$.  The smaller $\epsilon(\eta,A)$ becomes the
greater these restrictions are. 
In \cite{ekl} the
lower bound on the multiplicity which can be imposed at $\eta$ comes from
the Riemann--Roch theorem.  We improve upon this here by considering the
above mentioned obstructions and this
translates into a better lower bound for $\epsilon(\eta,A)$.  
The main difficulty in establishing the lower bound
$\epsilon(\eta,A) \geq \frac{1}{d-1}$ in general is that the set of points
where the bound $\epsilon(\eta,A) \geq \frac{1}{d}$ from \cite{ekl} fails
may contain divisors once $d \geq 3$.  
The reader only interested in the method employed and not the details of
counting should skip to \S 2 after going through Lemma \ref{ml}.

Finally, we would like to point out the similarity between the counting
methods used here and those employed by Faltings and W\"ustholz \cite{fw} to
reprove and extend the Schmidt subspace theorem.  In particular, the 
measure theoretic aspect of \cite{fw} is very closely related to our 
counting of jets.  The asymptotic dimensions
$h^0\left(X,nA \ts m_\eta^{n\alpha}/m_\eta^{n\alpha+1}\right)$, as $n \ra
\infty$, can be used to define a measure $\mu$
on $[0,m(A)]$, with $m(A)$ defined in \S 2.  Here $\mu((a,b))$ would
measure the asymptotic
cost of raising the multiplicity at $\eta$ from $a$ to $b$.

\medskip
\medskip
\noindent
{\bf Notation and Conventions}

\begin{itemize}
\item All varieties considered will be defined over the complex numbers
  $\cc$.
\item A point $x$ of an irreducible variety $X$ will be called very general if
    $x$ belongs to the complement of countably many closed, proper 
    subvarieties.
\item If $x \in X$ then $m_X \subset \oo_X$ is the maximal ideal sheaf of
$x$.
\item Suppose $V \subset X$ is an irreducible subvariety and $s \in H^0(X,L)$
      then $\ord_V(s)$ is the order of vanishing of $s$ along $V$.
\item If $L$ is a line bundle $\base(L)$ denotes the stable base locus of
          $L$, that is
$$
\base(L) = \{x \in X: s(x) = 0 \,\,\,\mbox{for all} \,\,s \in H^0(X,nL)\,\,\,
    \mbox{for all} \,\, n > 0 \}.
$$
\item  If $s \in H^0(X,L)$ then $Z(s)\subset X$ denotes the zero scheme of 
the section $s$.
\item  If $\alpha \in \rr$ then $\lfloor \alpha \rfloor$ denotes the 
  round--down of $\alpha$, that is the largest integer less than or equal
to $\alpha$.  Similarly, $\lceil \alpha \rceil$ denotes the round--up of
$\alpha$, the smallest integer greater than or equal to $\alpha$.  
\end{itemize}

\section{The Three--fold Case}

Before proving Theorem \ref{t} we review the strategy of \cite{ekl}.
Suppose then that $X$ is a smooth projective variety of
dimension $d$ and $A$ an ample line bundle on $X$.  Furthermore, let
$\eta$ be a very general point of $X$.  Ein, K\"uchle, and Lazarsfeld study
the linear series
$$
\left|kA \ts m_\eta^{k\alpha} \right|
$$
for various values of $\alpha$.

Roughly speaking, the argument of \cite{ekl} goes as
follows.  Suppose that $\epsilon(\eta,A) < \frac{1}{d}$ and
let $C_\eta$ be a curve with
\rn
\begin{eqnarray}
\frac{A \cdot C}{\mult_\eta(C_\eta)} < \frac{1}{d}.
\label{h}
\end{eqnarray}
Moreover, we assume that $C_\eta$ is chosen from a flat family 
${\cal F} \subset X \times T$
defined over a
smooth affine variety $T$ of dimension $d$ with a quasi--finite 
map $\phi: T \ra X$
where
$$
\frac{A \cdot C_t}{\mult_{\phi(t)}(C_t)} \, < \, \frac{1}{d}, \,\,\,\,\forall
t \in T.
$$
Consider, then, for $k$ sufficiently divisible, the linear series
\rn
\begin{eqnarray}
\left|kA \ts m_\eta^{k\alpha} \right|, \,\,\,k\alpha \in \zz.
\label{fff1}
\end{eqnarray}
If $k\alpha > k\epsilon(\eta,A)$ then by
\ref{h} the curve $C_\eta$ is in the base locus of this linear series.  
Using the fact that $C_\eta$ moves in the family ${\cal F}$ in order to
``differentiate in the parameter direction'', 
Ein, K\"uchle, and Lazarsfeld
then show that any divisor $D \in \left|kA \ts m_\eta^{k\alpha} 
\right|$ vanishes along $C_\eta$ to order at least $k\alpha -
k\epsilon(\eta,A)$.  
In particular, taking $\alpha > 2\epsilon(\eta,A)$ in \ref{fff1}, we see that
if $D \in \left|kA \ts m_\eta^{k\alpha} \right|$ then $D$ vanishes to order
greater than $k\epsilon(\eta,A)$ along $C_\eta$ and hence vanishes along
all curves in $C_t \in {\cal F}$ with $\phi(t) \in C_\eta$.  
The next step in the
argument is to show that a subfamily of
the curves $\{C_t\}_{\phi(t) \in C_\eta}$, defined over a constructible
subset $W \subset \phi^{-1}(C_\eta)$,
sweep out an irreducible surface $S_\eta \subset X$.
The argument is then iterated and the base locus of
$$
\left|kA \ts m_\eta^{k\alpha} \right|, \,\,\, \alpha > r\epsilon(\eta,A)
$$
is shown to contain an irreducible 
 subvariety of dimension $r$, swept out by an
$r-1$ dimensional subfamily of ${\cal F}$.  After iterating this 
argument $d$ times, a contradiction is reached because the linear series
$$
\left|kA \ts m_\eta^{kd\epsilon(\eta,A)+1} \right|
$$
is forced to be empty but by hypothesis $d\epsilon(\eta,A) < 1$ and
a simple argument using the Riemann--Roch theorem yields the
desired contradiction.  

\medskip
Fundamental to the argument of \cite{ekl} is  the following
``differentiation'' result:
\next
\begin{lemma}[\cite{ekl} Proposition 2.3]  
Suppose $\eta \in X$ is a very general point and let $W 
\subset X$ be an irreducible subvariety.  Let
$\pi: Y \ra X$ be the blow--up of $\eta$ with exceptional divisor $E$ and 
let $\tilde{W}$ is the strict transform of $W$ in $Y$.
Write
$$
\alpha(W) = \inf_{\beta \in \qq}\{\tilde{W} 
\subset \base(\pi^\ast A(-\beta E))\}.
$$ 
Suppose $\gamma > \alpha(W)$ and 
$0 \neq s \in H^0\left(X,nA \ts m_\eta^{n\gamma}\right)$.  Then
$$
\ord_W(s) \geq n\gamma - \lfloor \alpha(W)n + 1 \rfloor.
$$
\label{ml}
\end{lemma} 

\noindent
We have stated Lemma \ref{ml} in the form in which it will be used.  To 
obtain this version from \cite{ekl} Proposition 2.3, let $\Gamma \subset
 X \times T$ be the graph of $\phi: T \ra X$.
Let $p_1: X \times T \ra X$ and $p_2: X \times T \ra T$ denote the 
projections to the first and second factors respectively.  Consider
$$
\base \left(p_1^\ast(kA)\ts \ii_\Gamma^{\lfloor k\alpha(W) + 1\rfloor}\right).
$$
By hypothesis
for all $k > 0$ these subschemes contain an irreducible component $Z_k \subset
X \times T$
so that $Z_k \cap \pi_2^{-1}(t) \supset W$ for any $t$ with $\phi(t) = \eta$.
As $k \ra \infty$ one obtains a fixed subscheme $Z \subset X \times T$ with
$W \subset Z_t$ where $Z_t = Z \cap \pi_2^{-1}(t)$. 
According to \cite{ekl} Proposition 2.3, any section 
$\sigma \in H^0\left(X,p_1^\ast(kA) \ts \ii_\Gamma^{k\gamma}\right)$
must vanish to order at least $k\gamma - \lfloor k\alpha(W)+1\rfloor$ 
along $Z$.  Indeed, if not, after differentiating 
$k\gamma -\lfloor \alpha(W) + 1 \rfloor$ times we obtain $\sigma^\prime \in 
H^0\left(X,p_1^\ast(kA) \ts \ii_\Gamma^{\lfloor k\alpha(W) +1 \rfloor}\right)$
not vanishing along $Z$ and this is a contradiction.
Lemma \ref{ml} is the translation of this statement for the family $Z$ to 
the fibres $Z_t$.
Note that when applying Lemma \ref{ml} we often will assume, for
simplicity, that 
$$
\ord_W(s) \geq n(\gamma - \alpha(W)):
$$
indeed, for the asymptotic estimate on jets, the round--down and 1 are 
irrelevant.

\medskip

\medskip

\noindent
{\bf Proof of Theorem \ref{t}}  
Suppose that Theorem \ref{t} were false and
$$
\epsilon(\eta,A) < \frac{1}{2}.
$$
Thus through a very general point $\eta \in X$ there is a curve $C_\eta$
with $A\cdot C_\eta/\mult_\eta(C_\eta) < \frac{1}{2}$.  Choosing a suitable
family of such curves ${\cal F} \subset X \times T$ as above,
we claim that there is an open set $U \subset X$ such
that for each $x \in U$ there is an irreducible curve $C_x$ satisfying 
$A \cdot C_x = p$ and $\mult_\eta(C_x) = q$ with $p/q < 1/2$.  
If $\epsilon(\eta,A) = p/q$ and there is a curve $C_\eta$ through $\eta$
with $\mult_\eta(C_\eta) = q$ and $A\cdot C_\eta = p$ then this is satisfied.
If there were a  Seshadri exceptional surface $S$ at $\eta$, that is a surface
with $\frac{\deg_A(S)}{\mult_\eta(S)} = \epsilon(\eta,A)$, then
an immediate contradiction is obtained using Lemma \ref{ml}: choose 
$2p/q < \gamma < 1$ and $0 \neq \sigma \in H^0\left(X,nA \ts 
m_\eta^{n\gamma}/m_\eta^{n\gamma +1}\right)$.  According to Lemma \ref{ml}
$\ord_{S}(\sigma) \geq n\gamma - \lfloor pn/q + 1\rfloor$.  Since 
$\mult_\eta(S) \geq 2$ this is not possible for $n \gg 0$.  
Since there must either be
a Seshadri exceptional curve or a Seshadri exceptional surface when 
$\epsilon(\eta,A) < 1$ we must have a Seshadri excptional curve $C_\eta$
through $\eta$ as desired.

The goal of the proof is to estimate
\rn
\begin{eqnarray}
\lim_{n \ra \infty} 
\frac{h^0\left(X, nA \ts m_\eta^{\frac{3pn}{q}}\right)}{n^3}.
\label{y0}
\end{eqnarray}
We will show that this limit is positive and 
then we have a contradiction from \cite{ekl} which shows that the linear
series $\left|nA \ts m_\eta^{\frac{(3p+\alpha)n}{q}}\right|$ is empty for
any $\alpha > 0$.  
Let $\pi: Y \ra X$ be the blow--up of 
$X$ at $\eta$ with exceptional divisor $E$.  Choose a rational number $\alpha$
and a large positive integer $n$ with $n\alpha \in \zz$.  Then we have
\rn
\begin{eqnarray}
h^0(X,nA) &-& 
h^0\left(X,nA \ts m_\eta^{\alpha n} \right)  \nonumber \\
 & = & 
\sum_{k=0}^{\alpha n - 1}
\left(h^0\left(X,nA \ts m_\eta^{k} \right) -
 h^0\left(X,nA \ts m_\eta^{k+1} \right) \right)   \nonumber \\
& = &
\sum_{k=0}^{\alpha n - 1}
\left(h^0\left(Y,\pi^\ast(nA)(-kE)\right)  -
h^0\left(Y,\pi^\ast(nA)(-(k+1)E) \right) \right).
\label{ff1}
\end{eqnarray}
We have $E \simeq \pp^{2}$ and using \ref{ff1} and the exact sequence
\begin{eqnarray*}
0 \ra H^0\left(Y,\pi^\ast(nA)(-(k+1)E) \right)   \ra
H^0\left(Y,\pi^\ast(nA)(-kE)\right)  
\ra  H^0\left(E,\pi^\ast(nA)(-kE)\right)  
\end{eqnarray*}
we find
\rn
\begin{eqnarray}
h^0(X,nA) - 
h^0\left(X,nA \ts m_\eta^{\alpha n} \right)   =
\sum_{k=0}^{\alpha n - 1} h^0_Y\left(\pp^{2}, \oo(k)\right)
\label{ff2}
\end{eqnarray}
where $h^0_Y\left(\pp^{2}, \oo(k)\right)$ denotes the dimension of 
the subspace of $H^0\left(\pp^{2}, \oo(k)\right)$
coming via restriction
from $H^0\left(Y,\pi^\ast(nA)(-kE)\right)$.  
Our goal, then, is for each value of $k$, to bound
$h^0_Y\left(\pp^2,\oo(k)\right)$ from above.

\medskip

We next define critical numbers where
the base locus of $\left|kA \ts m_\eta^{k\alpha}\right|$ 
is forced to jump for numerical
reasons:
\begin{eqnarray*}
\alpha_1 = \frac{p}{q},\\
\alpha_3 = \frac{2p}{q}.
\end{eqnarray*}
There is also a more subtle jumping value between $\alpha_1$ and $\alpha_3$,
at least for $q$ sufficiently large, 
for which we require an extra definition.
Let $Z \subset \pp(T_\eta(X)) = \pp^2$ denote the zero--dimensional
subscheme of degree $q$ given by
$T_\eta(C_\eta)$.  Then one can define a Seshadri constant associated to 
$Z$ as follows.  Suppose $\psi: Y \ra \pp^2$ is a
birational map with $Y$ smooth and $\psi^{-1}(\ii_Z) = \oo_Y(-E)$.  Then 
$$
\epsilon(Z,\oo(1)) = \sup_r\{r\in \qq^+: \psi^\ast(\oo(1))(-rE)\,\,
\mbox{is nef}\}.
$$
Then, as we will see below,
 the base locus of $\left|kA \ts m_\eta^{k\alpha}\right|$ 
is forced to contain a
surface as soon as $\alpha > \alpha_2$ where $\alpha_2$ satisfies
$$
\frac{\alpha_2 - p/q}{2\alpha_2} = 
\epsilon \left(T_\eta(C_\eta),\oo_{\pp(T_\eta(X))}(1)\right).
$$
Note that a surface could enter the base locus of
$\left|kA \ts m_\eta^{k\alpha}\right|$  for $\alpha < \alpha_2$.  The numbers
we chose are the ``worst case scenario,'' the case where the linear
series $|kA|$ generates the {\em most possible} jets at $\eta$.  
If it generates
fewer jets, the numbers in the argument only improve.
We note here that the reader not interested in the counting details can
skip the analysis involving $\alpha_2$.  Indeed, when we prove Theorem
\ref{t} below there is enough room in the estimates so that the key result is
Lemma \ref{seshmult} which applies to the jet analysis once we have exceded
$\alpha_3$.  We included a more complete analysis both in order to reveal
the subtleties involved in counting and because in other cases the more
detailed analysis may be required.



By Lemma \ref{ml} we know that any
section of $\left|kA \ts m_\eta^{k\beta}\right|$ for $\beta > \alpha_1$ 
must vanish along $C_\eta$ to multiplicity at least 
$k(\beta - \alpha_1)$.  Once $\beta > \alpha_2$ we claim that
the base locus of $\left|kA \ts m_\eta^{k\beta}\right|$ must contain a 
surface $S$ which passes through $\eta$.  Indeed, if not, then choose
$s_1,s_2 \in \left|kA \ts m_\eta^{k\beta}\right|$ 
so that $T_\eta(Z(s_1))$ and $T_\eta(Z(s_1))$ meet properly inside
$T_\eta(X)$.
By Lemma \ref{ml}
we have $\mult_{C_\eta}(s_1) \geq k(\beta - \alpha_1)$  and
$\mult_{C_\eta}(s_2) \geq k(\beta - \alpha_1)$.
By Theorem A of \cite{els}, we have $f_1,f_2 \in \ii_C^{\lfloor k(\beta
- \alpha)/2 \rfloor}$ where $\ii_C$ is the ideal sheaf of $C$ and $f_1$
and $f_2$ are local equations for $s_1$ and $s_2$.  
Let $\pi: Y \ra X$ be the blow--up of $X$ at $\eta$ with exceptional
divisor $E$.  Let $D_1 = \pi^\ast(Z(s_1))(-k\beta E)|E$  and
$D_2 = \pi^\ast(Z(s_2))(-k\beta E)|E$.  We have $E \simeq \pp^2$ and
$D_1,D_2$ are curves of degree $k\beta$ meeting properly along 
$T_\eta(C_\eta)$, each with 
multiplicity at least $\lfloor k(\beta - \alpha_1)/2 \rfloor$ 
along $T_\eta(C_\eta)$.
Considering the pencil of divisors spanned by $D_1$ and $D_2$ shows that 
$\psi^\ast(\oo(k\beta))(-\lfloor k(\beta - \alpha_1)/2 \rfloor E)$ 
is nef where 
$\psi: Y \ra \pp^2$ is a resolution of $\ii_Z$ as above.  It follows that
$$
\epsilon(Z,\oo(1)) \geq \frac{\lfloor k(\beta -\alpha_1)/2 \rfloor}{k\beta} > 
\frac{\alpha_2 - p/q}{2\alpha_2},
$$
contradicting the definition of $\alpha_2$.  
Using Lemma \ref{ml} again, we conclude that there is a surface $S \subset X$
such that for $\beta > \alpha_2$ any divisor 
$D \in \left|kA \ts m_\eta^{k\beta}\right|$ must vanish along 
$S$ to order
 at least $k(\beta - \alpha_2)$.  Finally, let $S_\eta$ be the
surface swept out by $\{C_x\}_{x \in Z}$ with $Z \subset \phi^{-1}(C_\eta)$
the constructible subset considered above. By Lemma \ref{ml}
any divisor
 $D \in \left|kA \ts m_\eta^{k\beta}\right|$  must vanish along $S_\eta$ to
order at least $k(\beta - \frac{2p}{q})$.

\medskip

We are now prepared to bound
$h^0_Y(\pp^2, \oo(k))$ from above, using the information about the
order of vanishing of each section of $H^0_Y(\pp^2, \oo(k))$ along
$T_\eta(C_\eta)$, $T_\eta(S)$, and $T_\eta(S_\eta)$.
We will divide the estimate up into four cases
\begin{eqnarray*}
0 \leq k \leq n\alpha_1, \\
n\alpha_1  < k \leq n\alpha_2, \\
n\alpha_2  <k \leq  n\alpha_3,\\
n\alpha_3 < k \leq \frac{3np}{q}.
\end{eqnarray*}
We assume for simplicity that $\alpha_i \in \qq$ and $n\alpha_i \in \zz$.  
For those $\alpha_i$ which are irrational, it suffices in the argument 
below to replace $n\alpha_i$ by $\lfloor n\alpha \rfloor$. 
Note also that if $\alpha_2 > \alpha_3$, one simply eliminates the third
interval, replacing $\alpha_2$ by $\alpha_3$ in the second interval.  

For small values of $k$
 one expects $|nA|$ to generate all $k$--jets and
the estimate is
\rn
\begin{eqnarray}
h^0_Y\left(\pp^2,\oo(k)\right)
\leq 
\left( \begin{array}{c}    k+2 \\ 2 \end{array} \right), \,\,\,
0 \leq k \leq n\alpha_1.
\label{est1}
\end{eqnarray}
Next, for $n\alpha_1  < k \leq n\alpha_2$ 
any section $\sigma \in H^0_Y\left(\pp^2,\oo(k)\right)$ 
vanishes to order at least
$\lfloor(k-n\alpha_1)/2\rfloor$ along $T_\eta(C_\eta) \subset \pp(T_\eta(X))$
giving the estimate
\rn
\begin{eqnarray}
h^0_Y\left(\pp^2,\oo(k)\right)
\leq 
\left( \begin{array}{c}   k+2\\ 2 \end{array} \right) 
-q\left( \begin{array}{c}   \lfloor (k-n\alpha_1)/2 \rfloor + 1 \\ 
2 \end{array} \right) +
o(k^2), \,\,\,
n\alpha_1  < k \leq n\alpha_2.
\label{est2}
\end{eqnarray}
This is established in Lemma \ref{nonsep} below. 

Next suppose $n\alpha_2  < k \leq n\alpha_3$.  Let $\sigma \in 
H^0\left(X,nA \ts m_\eta^{k}\right)$.  For $n$ suitably divisible, write
$$
Z(\sigma) = aS + S^\prime, \,\,\,\mbox{with}\,\,\,\mult_\eta(S^\prime) = 
n\alpha_2:
$$
this is possible since $\sigma$ must vanish to order at least $k-n\alpha_2$
along the surface $S$.  Let 
$$ 
\rho: H^0\left(X,nA(-aS) \ts 
m_\eta^{n\alpha_2}\right) \ra H^0(\pp^2,\oo(n\alpha_2))
$$ 
be the restriction homomorphism.
Then by definition we have
$$
h^0_Y\left(\pp^2,\oo(k)\right) =  \dim({\rm Image}(\rho)).
$$
Using the construction in \cite{ekl} 3.8 we see that there exists an 
irreducible subvariety
$V \subset X \times T$ such that $S = V \cap
(X \times {t})$ for some $t$ with $\phi(t) = \eta$.  In particular, since
$\eta$ is a very general point it follows that there is a surface $S^\prime$
algebraically equivalent to $S$, not containing $\eta$, namely 
$S^\prime = V \cap (X \times {\xi})$ for a general point $\xi \in T$.  
Choose $r$ sufficiently large so that $rA + b(S-S^\prime)$ 
is very ample for all $b > 0$.
Choose $D \in |rA+a(S-S^\prime)|$ so that $D$ does not contain $\eta$ and
let $E = D+aS^\prime$.  Then tensoring by $E$ gives an injection
$$
\rho_E: H^0\left(X,nA(-aS) \ts 
m_\eta^{n\alpha_2}\right) \ra H^0\left(X,(n+r)A \ts  m_\eta^{n\alpha_2}
\right)
$$
which preserves multiplicity at $\eta$.  
We conclude that
\rn
\begin{eqnarray}
h^0_Y\left(\pp^2,\oo(k)\right) \leq h^0\left(X,(n+r)A \ts m_\eta^{n\alpha_2}/
m_\eta^{n\alpha_2+1}\right), \,\,\,n\alpha_2 < k \leq n\alpha_3.
\label{est3}
\end{eqnarray}

Finally suppose $n\alpha_3 < k \leq \frac{3pn}{q}$.  Suppose that
$\mult_\eta(\sigma) = k$, $\sigma \in H^0\left(X,nA\right)$.  
We know from Lemma \ref{ml}
that $\mult_{S_\eta}(\sigma) \geq k - n\alpha_3$.  
Since, according to Lemma \ref{seshmult} below 
$\mult_{C_\eta}(S_\eta) \geq 3$, we can write
$$
Z(\sigma) = aS_\eta + S^\prime
$$
with $\mult_\eta(S^\prime) = k - 3(k -n\alpha_3)$.  Arguing as in
the previous case then gives
\rn
\begin{eqnarray}
h^0_Y\left(\pp^2,\oo(k)\right) \leq h^0\left(X,(n+r)A \ts 
m_\eta^{3n\alpha_3 - 2k}/
m_\eta^{3n\alpha_3-2k+1}\right), \,\,\,n\alpha_3 < k \leq \frac{3pn}{q}.
\label{est4}
\end{eqnarray}

\medskip
We are now prepared to evaluate the limit \ref{y0} using \ref{ff2}.  
We assume to begin with that $q \geq 5$.  In particular this means that
$\epsilon \left(T_\eta(C_\eta), \oo_{\pp(T_\eta(X))}(1)\right) < 1/2$
and this guarantees that  $\alpha_2 < \alpha_3$.  
We divide the sum into four ranges of $k$ determined by our critical numbers
$\alpha_1,\alpha_2,\alpha_3$.  First, by \ref{est1}
$$
\lim_{n \ra \infty}
\frac{\displaystyle \sum_{k=0}^{n\alpha_1} h^0_Y(\pp^2,\oo(k))}{n^3}  \leq 
\frac{\alpha_1^3}{6}.
$$
Next, using \ref{est2} we see that 
$$
\lim_{n \ra \infty}
\frac{\displaystyle 
\sum_{k=n\alpha_1+1}^{n\alpha_2} h^0_Y(\pp^2,\oo(k))}{n^3}  \leq 
\frac{\alpha_2^3-\alpha_1^3}{6} - \frac{q(\alpha_2 -\alpha_1)^3}{24}.
$$
Using \ref{est2} and \ref{est3} gives
\begin{eqnarray*}
\lim_{n \ra \infty}
\frac{\displaystyle
\sum_{k=n\alpha_2+1}^{n\alpha_3} h^0_Y(\pp^2,\oo(k))}{n^3} & \leq &
\lim_{n \ra \infty}
\frac{\displaystyle \sum_{k=n\alpha_2+1}^{n\alpha_3} h^0\left(X,(n+r)A \ts 
m_\eta^{(n+r)(\frac{n\alpha_2}{n+r})}/m_\eta^{(n+r)(\frac{n\alpha_2}{n+r})+1}
\right)}{n^3} \\
& \leq & (\alpha_3-\alpha_2)\left(\frac{\alpha_2^2}{2} - 
\frac{q(\alpha_2-\alpha_1)^2}{8}\right).
\end{eqnarray*}
Finally, for $\sum_{k=n\alpha_3+1}^{3np/q} h^0_Y(\pp^2,\oo(k))$, using
\ref{est4} and arguing
as in the above case we can remove the $r$ which is fixed.
But the sum
$$
\sum_{k=n\alpha_3+1}^{3np/q} 
h^0\left(X,nA \ts 
m_\eta^{3n\alpha_3 - 2k}/
m_\eta^{3n\alpha_3-2k+1}\right)
$$
is simply every other term of the sum 
$\sum_{k=0}^{n\alpha_3}h^0_Y(\pp^2,\oo(k))$ and since
our upper bound for $h^0_Y(\pp^2,\oo(k))$ varies as a piecewise polynomial
this gives
$$
\lim_{n\ra \infty} \frac{\displaystyle \sum_{k=n\alpha_3+1}^{3np/q} 
h^0\left(\pp^2,\oo(k)\right)}{n^3}
\leq \frac{1}{2}\left(\frac{\alpha_2^3}{6} - 
\frac{q(\alpha_2 -\alpha_1)^3}{24}
+(\alpha_3-\alpha_2)\left(\frac{\alpha_2^2}{2} - 
q\frac{(\alpha_2-\alpha_1)^2}{8}\right)
\right).
$$
Combining all of the above estimates gives
\rn
\begin{eqnarray}
\lim_{n\ra \infty} \frac{\displaystyle \sum_{k=0}^{3np/q} 
h^0\left(\pp^2,\oo(k)\right)}{n^3}
&\leq& \frac{3}{2}\left(\frac{\alpha_2^3}{6} - 
\frac{q(\alpha_2 -\alpha_1)^3}{24}
+(\alpha_3-\alpha_2)\left(\frac{\alpha_2^2}{2} - 
q\frac{(\alpha_2-\alpha_1)^2}{8}\right)
\right)  \nonumber \\
&\leq&
\frac{3}{2}\left(\frac{\alpha_2^3}{6} 
+(\alpha_3-\alpha_2)\left(\frac{\alpha_2^2}{2}\right)\right).
\label{best}
\end{eqnarray}
Note that in the last inequality we have omitted two of the negative or
defect terms which were obtained by the detailed analysis above.  The reason
for this is that for $q$ sufficiently large, the estimate \ref{best} turns
out to be sufficient to establish Theorem \ref{t} while small values of
$q$ can be dealt with by hand.  
We included all of the counting
details nonetheless as this is the technical heart of this paper.

In order to compute the upper bound in \ref{best}, 
we need to know the value of $\alpha_2$. The Seshadri
constant $\epsilon\left(T_\eta(C_\eta),\oo_{\pp(T_\eta(X))}(1)\right)$ 
is, however, very difficult to compute and
so we look at the worst case scenario.  In particular, the bound in 
\ref{est2} increases until $x = \frac{np}{q-4} + O(1)$ and then 
decreases. The bound in \ref{est3} then repeats the last value for the
bound in \ref{est2} and then when one reaches \ref{est4} the values 
start to decrease.
 The $O(1)$ term will have no effect on the asymptotic estimate and
thus the worst case to cosider is $\alpha_2 = \frac{np}{q-4}$.  
With this value of $\alpha_2$ we need to assume that $q \geq 9$ in order
 to guarantee that $\alpha_2 < \alpha_3$.
We find then, using \ref{best}, 
\begin{eqnarray*}
\lim_{n\ra \infty} \frac{\displaystyle \sum_{k=0}^{3np/q} 
h^0\left(\pp^2,\oo(k)\right)}{n^3}
&\leq &\frac{3}{2}\left(\frac{p^3}{6(q-4)^3} 
+ \frac{p(q-8)}{q(q-4)}\left(
\frac{p^2}{2(q-4)^2} \right)\right) \\
&=& \frac{1}{6}\left(\frac{3p^3}{2(q-4)^2} + 
\frac{9p^3(q-8)}{2q(q-4)^3}\right)
\end{eqnarray*}
One checks that when $q \geq 10$ and $p/q < 1/2$ then
$$
\frac{1}{6}\left(\frac{3p^3}{2(q-4)^2} + 
\frac{9p^3(q-8)}{2q(q-4)^3}\right) < \frac{1}{6}.
$$
It follows from \ref{ff2} that when $q \geq 10$
$$
\lim_{n\ra \infty} \frac{h^0\left(X,nA \ts 
m_\eta^{\frac{3pn}{q}}\right)}{n^3} > 0
$$
and this concludes the proof of Theorem \ref{t} when $q \geq 10$.

If $q < 10$ then there are only four possibilities which are not eliminated
by \cite{ekl}, namely $p/q = 2/5, p/q = 3/7, p/q = 3/8,$ and $p/q = 4/9$.  
We outline
here how to eliminate the cases $p/q = 3/7$ and $p/q = 4/9$
 which are the most difficult of the
four.  The counting here goes as follows.  For $0 \leq k \leq np/q$ we
use \ref{est1}.  For $np/q < k \leq 2np/q$, we use the estimate in
\ref{f3} below.  Finally for $2np/q < k \leq 3np/q$, we use \ref{est4}.
This gives, in the case where $p/q = 3/7$,
$$
\lim_{n\ra \infty} \frac{\displaystyle \sum_{k=0}^{3np/q} 
h^0\left(\pp^2,\oo(k)\right)}{n^3}
\leq \frac{1}{6}\left(\frac{567}{686}\right).
$$
For $p/q = 4/9$ we find
$$
\lim_{n\ra \infty} \frac{\displaystyle \sum_{k=0}^{3np/q} 
h^0\left(\pp^2,\oo(k)\right)}{n^3}
\leq \frac{1}{6}\left(\frac{224}{243}\right).
$$

\medskip

\begin{lemma}
Suppose $C_\eta$ satisfies 
$$
\frac{A \cdot C_\eta}{\mult_\eta(C_\eta)} = \frac{p}{q} < \frac{1}{2}.
$$
Let $S_\eta$ be the surface swept out by $\{C_x\}_{x \in Z \subset
\phi^{-1}(C_\eta)}$.  Then
$$
\mult_{C_\eta}(S_\eta) \geq 3.
$$
\label{seshmult}
\end{lemma}

\noindent
{\bf Proof of Lemma \ref{seshmult}}  To see why $S_\eta$ must be
singular along $C_\eta$ note that for a general point $\xi \in \phi(Z) \subset
C_\eta$ 
there is a curve $C_\xi \subset S_\eta$ such that 
$$
\frac{A \cdot C_\xi}{\mult_\eta(C_\xi)} = \frac{p}{q} < \frac{1}{2}.
$$
If $S_\eta$ were smooth at a general point $\xi \in C_\eta$ then it would
follow that $\epsilon(\xi,A|S_\eta) < \frac{1}{2}$ and this is impossible
since Ein and Lazarsfeld \cite{el} have established that on a smooth
surface the set
of points where the
Seshadri constant can be less than one is at most countable.  In order
to refine this argument, let $\pi: X^\prime \ra X$ be an embedded resolution of
$S_\eta$.  For $\xi \in C_\eta$ general let
$\tilde{C}_\xi$ be the strict transform of $C_\xi$ in $X^\prime$ and write
$$
\pi^{-1}(\xi) \cap \tilde{C}_\xi = \{x_1,\ldots,x_r\}.
$$
Suppose moreover that $\psi: C \ra \tilde{C}_\xi$ is a desingularization with
$\psi^{-1}(\pi^{-1}(\xi)) = \left\{ y_1,\ldots,y_s \right\}$.  
Choose a linear series $|D|$ on $X$ with $D$ sufficiently positive so that if
$E \in |D|$ is a general member through $\xi$ then 
$i(x_j,\tilde{C}_\xi \cdot \tilde{E}:X^\prime) = \mult_{x_j}(\tilde{C}_\xi)$
for $1 \leq j \leq r$:
this is possible by \cite{F} 12.4.5.  
Then 
by \cite{F} 12.4.5  and 7.1.17 we have
$$
q = i(\xi,C_\xi \cdot D:X) = \sum_{j=1}^s 
{\rm ord}_{y_j}(\psi^\ast(\pi^\ast(D)) = \sum_{i=1}^r 
\mult_{x_i}\tilde{C}_\xi.
$$
Now we have $s = \mult_{C_\eta}(S_\eta)$ and thus if 
$s = 2$ we find that for $i = 1$ or $i = 2$
$$
\frac{\pi^\ast(A) \cdot \tilde{C}_\xi}{\mult_{x_i}(\tilde{C}_\xi)} < 1.
$$
Let $\tilde{S} \subset X^\prime$ be the resolution of $S_\eta$.  The
curves $\tilde{C}_\xi$ move in a
one parameter family along the surface $\tilde{S}$ and thus
$\{x \in \tilde{S}: \epsilon(x,\pi^\ast(A)) < 1\}$ is not countable, 
violating the main result of \cite{el}.  Note that in \cite{el} Ein and 
Lazarsfeld state the main result for ample line bundles but the proof holds
unchanged for a big and nef line bundle. Indeed, the only point where
\cite{el} uses ampleness is to show that curves with bounded degree relative
to the appropriate ample bundle $A$ 
move in finitely many families but this also holds
more generally when $A$ is big and nef.

\medskip
\begin{lemma}
Suppose $\mult_\eta(C_\eta) = q$ and $n\alpha_1 <k \leq n\alpha_2$.  Then
$$
h^0_Y\left(\pp^2,\oo(k)\right)
\leq 
\left( \begin{array}{c}    k+2\\ 2 \end{array} \right) 
-q\left( \begin{array}{c}    \lfloor (k-n\alpha_1 + 1)/2\rfloor \\ 
2 \end{array} \right) + o(k^2), \,\,\,
n\alpha_1  < k \leq n\alpha_2.
$$
\label{nonsep}
\end{lemma}

\noindent
{\bf Proof of Lemma \ref{nonsep}}
Let $A = \oo_{\pp^2}(1)$ where $\pp^2 = \pp(T_\eta(X))$ and let
$Z \subset \pp^2$ be the projectivized tangent cone of $C_\eta$ at $\eta$.  
By definition of $\epsilon(Z,A)$, given $\delta > 0$ so that
$\epsilon(Z,A) - \delta \in \qq$ for all $n > 0$ sufficiently large and 
divisible the 
evaluation map
$$
H^0(\pp^2,\oo(n)) \ra 
H^0\left(\pp^2,\oo(n) \ts \oo_{\pp^2}/\ii_Z^{n(\epsilon(Z,A) - \delta)}\right)
$$
is surjective.
According to \cite{F} Example 4.3.4
$$
\ell(\oo_X/\ii_Z^{r}) = \frac{qr^2}{2} + O(r).
$$
In particular, for $n\alpha_1 < k \leq n\alpha_2$ we see that that
Lemma \ref{nonsep} holds since any section of $H^0_Y\left(\pp^2,\oo(k)\right)$
vanishes to order at least $\lfloor(k - n\alpha_1)/2\rfloor$ along $Z$.  


\medskip
\medskip
\medskip

\section{Counting jets in the general case}

In addition to the computational complexity involved in estimating
$h^0\left(X,nA \ts m_\eta^{k}\right)$ for different values of $k$, the
central difficulty in the higher dimensional case is Lemma \ref{seshmult} which
uses the fact that on a surface $X$ the set $\{x\in X:\epsilon(x,A) < 1\}$ is
countable
for an ample line bundle $A$.  In particular, it is critical for 
Lemma \ref{seshmult}
that this set contains no divisor and this is not known in higher dimension.
In order to prove Theorem \ref{main},
we begin by recalling a key definition from \cite{N}.  
\begin{definition}  
For an ample $\qq$--divisor $A$ on a smooth surface $X$ we let
$$
m(A) = \sup_{D \equiv A}\left\{\mult_\eta(D)\}\,|
\,\,D \in \,\,{\rm Div}(X) \ts
\qq \,\,\mbox{effective}\right\}:
$$
here $\equiv$ denotes numerical equivalence.
\label{d2}
\end{definition}

\noindent
The importance of Definition \ref{d2} lies in the following simple result:

\begin{lemma}  Suppose $X$ is a projective variety of dimension $d$ 
and $A$ an ample line bundle on $X$.  If 
$$
\epsilon(\eta,A) < \frac{m(A)}{d}
$$
then given $\delta > 0$
there exists an irreducible proper subvariety $Y \subset X$, of dimension
at least one, such that 
$$
\epsilon(\xi,A|Y) <  \epsilon(\eta, A) + \delta:
$$
here $\xi$ is a very general point of $Y$ and $A|Y$ is the restriction of $A$
to $Y$.
\label{l0}
\end{lemma}

\noindent
{\bf Proof of Lemma \ref{l0}}  Suppose, to the contrary, that for some
$\delta > 0$
$$
\epsilon(\xi,A|Y) \geq \epsilon(\eta,A) + \delta
$$
for every irreducible $Y \subset X$.  As above, 
following \cite{ekl} 3.4, given $\epsilon > 0$ there is a 
a family of curves ${\cal F} \subset X \times T$ with
\rn
\begin{eqnarray}
\alpha = \frac{A \cdot C_t}{\mult_{\phi(t)}(C_t)} < \epsilon(\eta,A) + 
\epsilon,
\,\,\,t \in T.  
\label{tg}
\end{eqnarray}
This gives a chain of subvarieties
\rn
\begin{eqnarray}
V_1 \subset V_2 \subset \cdots \subset V_d
\label{l1}
\end{eqnarray}
where $V_1 = C_\eta$ for a very general point $\eta$ and $V_{i+1} = C(V_i)$
in the notation of \cite{ekl} Lemma 3.5.1: 
in particular, $V_{i+1}$ is obtained by 
adjoining curves in the family ${\cal F}$.  By
\cite{ekl} Lemma 3.5.1 each $V_i$ is irreducible. 
According to Lemma \ref{ml} any section 
$$
s \in H^0\left(X, nA \ts m_\eta^{2n\alpha + 2}\right)
$$
vanishes along $C_\eta$ to order at least $n\alpha + 1$ 
and hence vanishes along
$V_2$.  Proceding inductively using Lemma \ref{ml} we find that if
$s \in H^0\left(X, nA \ts m_\eta^{dn\alpha + d}\right)$ then
\rn
\begin{eqnarray}
s|V_d = 0.
\label{uv}
\end{eqnarray}
By hypothesis, $m(A) > d\epsilon(\eta,A)$ and thus, shrinking $\epsilon$ 
in \ref{tg} if
necessary, 
we can assume that $s$ is not indentically
zero in \ref{uv} and thus $\dim (V_d) \leq d-1$. It 
follows from \ref{l1} that for some $1 \leq r \leq d-1$ we must have
$$
V_r = V_{r+1}.  
$$
In particular for a general, hence smooth,  point $\xi \in V_r$, we find
a curve $C_\xi \subset V_r$ with
$$
\frac{\mult_\xi(C_\xi)}{A \cdot C_\xi} = \alpha.
$$
Hence
$$
\epsilon(\xi,A|V_r) \leq \alpha < \epsilon(\eta,A) + \delta.
$$
Thus we can take $Y = V_r$ and this proves Lemma \ref{l0}.

\medskip

We will derive Theorem \ref{main} from Lemma \ref{l0} and the following
result.  
\begin{lemma}
Suppose $X$ is a projective variety of dimension $d \geq 4$ 
and $A$ an ample line
bundle on $X$.  Then either
$$
m(A) > 1 + \frac{1}{3d}
$$
or 
$$
\epsilon(\eta,A) > \frac{1}{d} + \frac{1}{3d^2}.
$$
\label{l2}
\end{lemma}

\noindent
{\bf Proof of Lemma \ref{l2}}  
The proof
of Lemma \ref{l2} follows closely the method of \S 1 though the counting
is much simpler.  Let $\pi: Y \ra X$ be the blow--up of 
$X$ at $\eta$ with exceptional divisor $E \simeq \pp^{d-1}$.  
Choose a rational number $\alpha$,
and a large positive integer $n$
so that $n\alpha \in \zz$.  Then we have, with the same notation
 as above
\rn
\begin{eqnarray}
h^0(X,nA) - 
h^0\left(X,nA \ts m_\eta^{\alpha n} \right)   =
\sum_{k=0}^{\alpha n - 1} h^0_Y\left(\pp^{d-1}, \oo(k)\right)
\label{f2}
\end{eqnarray}
where, as above,
 $h^0_Y\left(\pp^{d-1}, \oo(k)\right)$ denotes the dimension of 
the subspace of $H^0\left(\pp^{d-1}, \oo(k)\right)$
coming via restriction
from $H^0\left(Y,\pi^\ast(nA)(-kE)\right)$.  

Suppose that $x \in \pp^{d-1} = \pp(T_\eta(X))$ is a tangent vector to
$C_\eta$ at $\eta$.  Then for $k > \epsilon(\eta,A)n$, Lemma \ref{ml} implies
\rn
\begin{eqnarray}
h^0_Y\left(\pp^{d-1}, \oo(k)\right) \leq
h^0\left(\pp^{d-1}, \oo(k) \ts m_x^{\lceil k-\epsilon(\eta,A)n -1\rceil}
\right):
\label{f3}
\end{eqnarray}
Combining \ref{f2} and \ref{f3} and taking the limit as $n \ra \infty$
we find
\begin{eqnarray}
\lim_{n \ra \infty} \frac{h^0(X,nA)- h^0(X,nA \ts m_\eta^{\alpha n})}{n^d}
&\leq& 
\frac{1}{(d-1)!}\int_0^\alpha \left(x^{d-1} - 
{\rm max}\left\{0,(x-\epsilon(\eta,A))^{d-1}\right\}\right)  \nonumber \\
&=& \frac{\alpha^d - (\alpha - \epsilon(\eta,A))^d}{d!}.
\rn
\label{f4}
\end{eqnarray}
According to \ref{f4}, if 
$$
\alpha^d - (\alpha - \epsilon(\eta,A))^d < 1
$$
then $\lim_{n \ra \infty} 
\frac{h^0\left(X,nA \ts m_\eta^{\alpha n}\right)}{n^{d}} > 0$ and
hence $m(A) > \alpha$. Suppose then that
$\epsilon(\eta,A) \leq \frac{1}{d}+ \frac{1}{3d^2}$ and let
$\alpha = 1 + \frac{1}{3d}$.  Then we find
$$
\lim_{d \ra \infty}\left( \left(\frac{3d+1}{3d}\right)^d - 
\left(\frac{3d+1}{3d} - \frac{3d+1}{3d^2}\right)^d\right) \leq e^{1/3} -
e^{-2/3} < 0.9.
$$
Thus we see that for all $d$ sufficiently large $m(A) > 1+ \frac{1}{3d}$ if 
$\epsilon(\eta,A) \leq \frac{1}{d}
+ \frac{1}{3d^2}$.  Elementary calculus suffices to show that
this also holds for all $d \geq 4$  and this establishes Lemma \ref{l2}.

\medskip

\noindent
{\bf Proof of Theorem \ref{main}} Suppose to the contrary that
\rn
\begin{eqnarray}
\epsilon(\eta,A) \leq \frac{1}{d} + \frac{1}{3d^2}.
\label{g1}
\end{eqnarray}
By  Lemma \ref{l2}  $m(A) > 1 + 1/3d$  and thus 
$$
\epsilon(\eta,A) < \frac{m(A)}{d}.
$$
By Lemma \ref{l1}
given $\delta > 0$  there
is a proper subvariety $Y \subset X$ such that for a very general point $\xi 
\in Y$
$$
\epsilon(\xi,A|Y) \leq \frac{1}{d}+ \frac{1}{3d^2} + \delta.
$$
But by the main theorem of \cite{ekl}, we know that $\epsilon(\xi,A|Y) \geq
\frac{1}{\dim(Y)} \geq \frac{1}{d-1}$ and this is a contradiction for 
$\delta$ sufficiently small.  

\medskip
Note that above in \ref{f3} we have not counted carefully: in particular,
the curve $C_\eta$ is singular at $\eta$ and thus the tangent
space to $C_\eta$ will be more than a single point with multiplicity one.  Thus
the counting can be improved considerably here but we were unable to
obtain a significant quantitative improvement in the final result by checking
this counting more carefully.

\begin{tabbing}
Department of Mathematics and Statistics\\
University of New Mexico\\
Albuquerque, New Mexico 87131\\
{\em Electronic mail:} nakamaye@math.unm.edu
\end{tabbing}

\end{document}